\documentclass[twoside, 12pt]{article}
\usepackage{filecontents}
\usepackage{amsfonts}
\usepackage{amsmath}
\usepackage{amssymb}
\usepackage{blkarray}
\usepackage{multirow}
\usepackage{graphics}
\usepackage[noadjust]{cite}
\usepackage{amsthm} 
\usepackage{varioref}
\usepackage[colorlinks=true,allcolors=blue]{hyperref}
\usepackage[nameinlink,capitalize,noabbrev]{cleveref}
\crefname{equation}{}{}

\usepackage{pgf,tikz}
\usepackage{mathrsfs}
\usetikzlibrary{arrows}
\definecolor{wrwrwr}{rgb}{0.3803921568627451,0.3803921568627451,0.3803921568627451}
\definecolor{rvwvcq}{rgb}{0.08235294117647059,0.396078431372549,0.7529411764705882}
\newcommand{\R}{\mathbb{R}}
\newtheorem{theorem}{Theorem}[section]
\newtheorem{remark}[theorem]{Remark}
\newtheorem{example}[theorem]{Example}
\newtheorem{lemma}[theorem]{Lemma}
\newtheorem{corollary}[theorem]{Corollary}
\newtheorem{definition}[theorem]{Definition}
\newtheorem{proposition}[theorem]{Proposition}

\newtheorem{conjecture}{Conjecture}

\newcommand{\bea}{\begin{eqnarray}}
\newcommand{\eea}{\end{eqnarray}}

\newcommand{\comment}[1]{}
\def \bpm{\begin{pmatrix}}
\def \epm{\end{pmatrix}}
\def \bd{\begin{definition}}
\def \ed{\end{definition}}
\def \bcc{\begin{conjecture}}
\def \ecc{\end{conjecture}}
\def \bt{\begin{theorem}}
\def \et{\end{theorem}}
\def \bl{\begin{lemma}}
\def \el{\end{lemma}}
\def \bc{\begin{corollary}}
\def \ec{\end{corollary}}
\def\be#1\ee{\begin{align}#1\end{align}}
\def\beq #1\eeq {\begin{align*}#1\end{align*}}
\def \ben{\begin{enumerate}}
\def \een{\end{enumerate}}
\def \ba{\begin{array}}
\def \ea{\end{array}}
\def \bp{\begin{proposition}}
\def \ep{\end{proposition}}
\def \bx{\begin{example}}
\def \ex{\end{example}}
\def \br{\begin{remark}}
\def \er{\end{remark}}
\def \bdsc{\begin{description}}
\def \edsc{\end{description}}
\def\pf{{\it \bf Proof. }}
\def \qed {\hfill \vrule height6pt width6pt depth0pt}
\def\hs{\hspace{.3cm}}
\def\vs{\vskip .3cm}
\def\ds{\displaystyle}
\def\1{1\!\!1}
\def\J{\mathbb{J}}
\def\x{\mathbf{x}}
\def\r{\rho}
\def\D{\mathcal{D}}

\title{On Pareto eigenvalue of distance matrix of a graph}
\author{Milan Nath,~Deepak Sarma$^{1,}$\footnote {Corresponding author \newline Email addresses: \url{milan@tezu.ernet.in} (M. Nath), \url{deepaks@tezu.ernet.in} (D. Sarma).\newline $^1$The financial assistance for this author was provided by CSIR, India, through JRF.}
 \\Department of Mathematical Sciences, \\ Tezpur University, Tezpur-784028, India.}

\date{}

\begin{document}
\maketitle

\vs

\begin{abstract}
In this article, we study Pareto eigenvalues of distance matrix of connected graphs and show that the non zero entries of every distance Pareto eigenvector of a tree forms a strictly convex function on the forest generated by the vertices corresponding to the non zero entries of the vector. Besides we find the minimum number of possible distance Pareto eigenvalue of a connected graph and establish lower bounds for $n$ largest distance Pareto eigenvalues of a connected graph of order $n.$ Finally, we discuss some bounds for the second largest distance Pareto eigenvalue and find graphs with optimal second largest distance Pareto eigenvalue.


\end{abstract}

\maketitle
\bigskip

\section{Introduction and terminology}

\hs All graphs considered here are finite, undirected, connected and simple.
Let $G$ be a  graph on vertices $1,2,\ldots,n$. At times, we use $V(G)$ and
$E(G)$ to denote the set of vertices and the set of edges of $G$, respectively. For $i, j\in V(G)$,
the {\em distance} between $i$ and $j$, denoted by $d_G(i,j)$ or simply $d_{ij}$, is the length of
a shortest path from $i$ to $j$ in $G$. The {\em distance matrix} of $G$, denoted by
$\D(G)$ is the $n\times n$ matrix with $(i,j)$-th entry $d_{ij}$.
For a column vector $x=(x_1,\ldots,x_n)^T\in\mathbb{R}^n$ we have
\beq x^T\D(G)x=\sum_{1\le i<j\le n}d_{ij}x_ix_j. \eeq

If vertices $i$ and $j$ are adjacent, we write $i\sim j.$ Edges $e_1, \,e_2$ in a graph $G$ are said to be incident if they have a common vertex.
If $V_1\subseteq V(G)$ and $E_1\subseteq E(G),$ then by $G-V_1$ and $G-E_1$ we mean the graphs obtained from $G$ by deleting the vertices in $V_1$ and the edges $E_1$ respectively. In particular case when $V_1=\{u\}$ or $E_1=\{e\},$ we simply write $G-V_1$ by $G-u$ and $G-E_1$ by $G-e$ respectively. $K_n-e$ is the graph obtained from $K_n$ by removing any one edge of it. Degree of a vertex $v$ in a graph $G$ will be denoted by $d_v.$ By a pendent vertex of a graph we mean a vertex of degree 1. For a connected graph $G$ a block $B$ is said to be a pendent block if exactly one vertex of $B$ is a cut vertex of $G.$ A \textit{quasipendent vertex}
is a vertex which is adjacent to a pendent vertex. The \textit{transmission}, denoted by $Tr(v)$ of a vertex $v$
is the sum of the distances from $v$ to all other vertices in $G$. The \textit{Wiener index}, denoted by $W$ of a connected graph $G$ is defined as $W={\ds \frac{1}{2}\sum_{v\in V(G)}}Tr(v).$  By $P_n, K_n$ and $C_n$ we mean the usual path graph, complete graph and cycle graph with $n$ vertices respectively. The \textit{diameter} of a connected graph $G$ denoted by diam$(G)$ is the maximum distance between any two vertices in $G,$ i.e. diam$(G)$ is the largest entry of $\D(G).$ A \textit{clique} of a graph is a maximal complete subgraph and \textit{clique number} of a graph is the order of a maximal clique. We denote \textit{clique number} of a graph $G$ by $\omega(G).$ For $1\leq p\leq \omega$ by $K_\omega ^p$ we denote a graph obtained by joining one vertex to $p$ vertices of $K_\omega.$ By $S_n$ we mean the usual star graph $K_{1,n-1}$ and by $S_n^+$ we represent the graph so that  $S_n^+-e=S_n.$ The graph obtained from $G$ and $H$ by identifying $u\in G$ and $v\in H$ is denoted by $G_u*H_v$ or simply by $G*H$ when there is no confusion of the vertices. We write $H_{u,v}$ to denote a graph of order $n$ with $u,v\in V(H_{u,v})$ so that $d_u=n-1$ and each vertex in $V(H)- \{u, v \}$ has same vertex degree and same transmission. 


By \textit{spectral radius} of a symmetric matrix $M$, we mean its largest eigenvalue and denote it by $\rho(M).$
Note that for a connected graph $G$, $\D(G)$ is irreducible nonnegative
matrix. Thus by the Perron-Frobenius theorem, $\rho(\D)$ is simple, and there is a
positive eigenvector of $\D(G)$ corresponding to $\rho(\D)$.
Such eigenvectors corresponding to $\rho(\D)$ is called {\em Perron vector}
of $\D(G)$. By an eigenvector we mean a unit eigenvector and by $\mathbb{M}_n,$ we denote the class of all real matrices of order n. We use the notation $A\geq 0$ to indicate that each component of the matrix $A$ is nonnegative. Furthermore in places we write $A\geq B$ to mean $A-B\geq 0.$

\bd
A real number $\lambda$ is said to be a Pareto eigenvalue of $A\in \mathbb{M}_n$ if there exists a nonzero vector $\x(\geq 0)\in \R^n$ such that \\
\beq A\x \geq\lambda \x \quad and \quad \lambda=\frac{\x^TA\x}{\x^T\x}, \eeq
also we call $\x$ to be a Pareto eigenvector of $A$ associated with Pareto eigenvalue $\lambda$.
\ed

Pareto eigenvalues are also known as complementarity eigenvalues. Fernandes {\it et al.} \cite{fjt17}  and Seeger \cite{see18} studied the Pareto eigenvalues of adjacency matrix of a graph.
\par We now outline the contents of this article. In \cref{sec2}, we introduce distance Pareto eigenvalue (eigenvector) of a connected graph and show that non zero entries of every distance Pareto eigenvector of a tree forms a strictly convex function on the forest generated by vertices corresponding to the non zero entries of the vector. We also find the complete distance Pareto spectrum for some special class of graphs like complete graph, Star graph etc. Partial distance Pareto spectrum of graphs with given diameter or given clique number are also supplied. Besides we find the minimum number of possible distance Pareto eigenvalue of a connected graph and establish lower bounds for $n$ largest distance Pareto eigenvalues of a connected graph of order $n$, equality conditions have also been established. In \cref{sec5}, we discuss the second largest distance Pareto eigenvalue of a connected graph, specially we give some bounds for it and find graphs with optimal second largest distance Pareto eigenvalue.

\section{Distance Pareto eigenvalue of a connected graph}\label{sec2}

\bd
Distance Pareto eigenvalue of a connected graph $G$ is a Pareto eigenvalue of the distance matrix of $G.$
\ed
\hs Multiplicity of Pareto eigenvalue of a matrix is not considered. We denote the $k^{th}$ largest and $k^{th}$ smallest distance Pareto eigenvalue of a connected graph $G$ by $\rho_k(G)$ and $\mu_k(G)$ respectively. We simply write them by $\rho_k$ and $\mu_k$ when the graph under consideration is understood from the context. Besides we use $\Pi(G)$ to denote the set of all distance Pareto eigenvalues of a connected graph $G.$
\par We write $[n]=\{1,2,\ldots ,n\}$ and for an $n\times n$ matrix $A$ and $S\subset [n],$ we reserve the symbol $A(S)$ for the principal submatrix of $A$ obtained by deleting rows and columns of $A$ corresponding to $S.$  In particular if $\D(G)$ is the distance matrix of a graph $G$ then by $\D(i)$ we will denote the principal submatrix of $\D$ obtained by deleting row and column corresponding to vertex $i$ of $G.$ By $\1$ we denote the column vector of all ones and by $\J$ the matrix of all ones of appropriate size.
We now recall some known results which will be used.

 \bl{ [Weyl's Inequalities]}\cite{hj05}\label{w} Let $\lambda_i(M)$ denote the $i$-th largest eigenvalue of a real symmetric matrix $M.$ If $A$ and $B$ are two real symmetric matrices of order $n,$
then $$\lambda_1(A)+\lambda_i(B)\ge \lambda_i(A+B)\ge \lambda_n(A)+\lambda_i(B) \mbox{ for } i=1,2,\ldots, n.$$\el

\bl \cite{min88}\label{mi} If $A$ is an irreducible matrix and $A\geq B\geq 0, A\neq B,$ then $\r(A)>\r(B).$
\el

\bl { [Cauchy's Inequalities]}\cite{hj90} \label{cauchy} Let $A=\bpm
B & y \\
y^T & a
\epm \in \mathbb{M}_{n+1}$ with $y\in \mathbb{C}^n, a \in \R$ and $B\in \mathbb{M}_n$ be symmetric. Then
\beq \lambda_1(A)\leq \lambda_1(B)\leq  \lambda_2(A)\leq \cdots \leq \lambda_n(A)\leq \lambda_n(B)\leq \lambda_{n+1}(A), \eeq
where $\lambda_i(N)$ is the $i^{th}$ smallest eigenvalue of $N.$ Also  $\lambda_i(A)=\lambda_i(B)$ if and only if there is a non zero vector $z \in \mathbb{C}^n$ such that $Bz=\lambda_i(B)z, y^Tz=0$ and $Bz=\lambda_i(A)z.$
\el

\bl \cite{hj05}\label{non} Let $A,B \in \mathbb{M}_n$ be nonnegative where $A$ is irreducible and $B$ is non zero, then $\r(A+B)>\r(A).$
\el

\bl\cite{hj05}\label{r} If $A$ is a symmetric $n\times n$ matrix with $\lambda_1$ as the largest eigenvalue then for any normalized vector
$\x\in \R^n (\x\ne 0)$,
\beq \x^t A \x\le \lambda_1.\eeq The equality holds if and only if $\x$ is an eigenvector of $A$ corresponding to the eigenvalue $\lambda_1.$\el

Putting $\x=\frac{1}{\sqrt{n}}\1$ in \cref{r} we get the following result as a Corollary.

\bc \label{avg} If $A$ is a symmetric $n\times n$ matrix, then
\be \label{av} \r(A)\ge \bar{R}, \ee \\
 where $\bar{R}$ is the average row sum of the matrix $A.$ The equality in \cref{av} holds if and only if all the row sums of $A$ are equal.
\ec

\bt \cite{sv11}\label{imp}
The scalar $\lambda \in \R$ is a Pareto eigenvalue of $A\in \mathbb{M}_n$ if and only if there exists a nonempty set $J\subset \{ 1, 2, \ldots, n\}$ and a vector $\xi\in \R^{|J|}$ such that

\beq
A^J\xi &=\lambda\xi, \\
\xi_j  &>0 \quad \text{ for all }j \in J,  \nonumber \\
 \sum_{j\in J}a_{i,j}\xi_j & \geq 0 \quad \text{ for all } i\notin J. \nonumber
\eeq
 Furthermore, a Pareto eigenvector $\x$ associated to $\lambda$ is constructed by setting
\be
\x_j =
  \begin{cases}
   \xi_j & \text{if } j\in J, \nonumber \\
   0       & otherwise.
  \end{cases}
\ee
\et

\par From \cref{imp}, we get the following result similar to that of \cite[Theorem 1]{see18}.
\bt \label{main} The distance Pareto eigenvalues of a connected graph $G$ are given by

$\Pi(G)=\{ \r(A): A\in M\},$
where $M$ is the class of all principal sub-matrices of $\D(G).$
\et

From \cref{main}, for any connected graph $G$ we see that every Pareto eigenvalue of any principal sub-matrix of $\D(G)$ is the distance Pareto eigenvalue of $G$. Therefore as a consequence we get that if $G$ is a connected graph and $H$ is a block of $G$, then every distance Pareto eigenvalue of $H$ is also a distance Pareto eigenvalue of $G.$ Also if $G$ is a connected graph and $H$ is a subgraph of $G$ obtained by removing one or more pendent blocks of $G,$ then $\lambda \in \Pi(H)$ implies $\lambda \in\Pi(G).$ Again if $\omega$ is the clique number of the graph $G$, then $\J_{\omega}-I_{\omega}$ is a principal sub-matrix of $\D(G)$ and  eigenvalues of $\J_{\omega}-I_{\omega}$ are $0, 1, \ldots , \omega-2$ and $\omega-1.$ Therefore $0,1,\ldots,\omega-1$ are distance Pareto eigenvalues of G. Similarly if there are $p$ vertices in a graph $G$ which are at a distance $k$ from each other, then $k(\J_p-I_p)$ is a principal sub-matrix of $\D(G)$ and therefore  eigenvalues of $k(\J_p-I_p)$ i.e. $0,k,\ldots,k(p-2)$ and $k(p-1)$ are distance Pareto eigenvalues of G. Besides from \cref{main}, we get the following lemma which states that for a connected graph the distance spectral radius and the largest distance Pareto eigenvalue coincide.

\bl \label{max} The largest distance Pareto eigenvalue of a connected graph is the distance spectral radius of the graph.
\el

\bd \cite{bkns17} Let $T$ be a tree with $V(T) = \{1, \ldots, n\}, n \ge 3,$ and let $f:V(T) \longrightarrow
[0, \infty).$
Then $f$ is said to be {\it convex} if for any distinct $i,j,k \in V(T)$ with $i \sim j, j \sim k$
we have $2f(j) \le f(i) + f(k)$ and {\it strictly convex} if  $2f(j) < f(i) + f(k).$
Also $f$ is said to be {\it quasiconvex} if for any distinct $i,j,k \in V(T)$ with $i \sim j, j \sim k$
we have $f(j) \le \max\{f(i), f(k)\}$ and {\it strictly quasiconvex} if  $f(j) < \max\{f(i) , f(k)\}.$ A function defined on the vertices of a forest $F$ is said to be convex, strict convex, quasiconvex or strictly quasiconvex if it is so in any subtree of $F.$ \ed

 \bl \cite{bkns17} \label{lem1}
 Let $T$ be a tree with $V(T) = \{1, \ldots, n\}$ and let $f:V(T) \longrightarrow
 [0,\infty)$ be either strictly convex or strictly quasiconvex. Then $f$ attains its minimum, either at a
unique vertex, or at two adjacent vertices. Furthermore,
\begin{description}
\item {(i)} if $f$ attains its minimum at the unique vertex $i,$ then for any path
$i=i_1-i_2-\cdots-i_k,$  starting at $i,$  $f(i_1) < f(i_2) < \cdots < f(i_k).$

 \item {(ii)} if $f$ attains its minimum at the two adjacent vertices
 $i$ and $j,$ then $f$ is strictly increasing along any path starting at $i,$
 and not containing $j,$ or starting at $j,$ and not containing $i.$

 \end{description}
\el

\bl \cite{bkns17} \label{cor11}
Let $\D$ be the distance matrix of a tree $T$ and let $f:[0,\infty) \longrightarrow [0,\infty)$ be a strictly increasing
   convex function. Then the
Perron vector of $f(\D)$ is strictly convex on $T$.
\el

\bt If $\x$ is a distance Pareto eigenvector of a tree $T$ and $V_x=\{v\in T: \x_v>0 \}$ then the non zero components of $\x$ form a strictly convex function on the subgraph of $T$ generated by $V_x.$
\et
\pf If $\x$ is a distance Pareto eigenvector of tree $T$ with corresponding distance Pareto eigenvalue $\lambda $ then
from \cref{main} and \cref{imp}, $\lambda$ is a spectral radius for some principal sub-matrix $M$ of $\D(T)$ and $x_v>0$ if and only if $v^{th}$ row(column) of $\D(T)$ is in $M.$ If $\lambda=0,$ then cardinality of $V_x$ is 1 and therefore the result is trivially true. So we assume $\lambda\neq 0.$ Now if the subgraph of $T$ generated by $V_x$ is again a tree, then the result follows from \cref{cor11}.
\par Otherwise we consider a subtree $T'$ of the subgraph (forest) of $T$ generated by $V_x.$ Let $u,v,w$ be vertices of $T'$ such that $u \sim v\sim w.$
Let $e_1 = \{u,v\},   e_2 =  \{v,w\}$ be edges of $T'.$  Let $T_u, T_v$ and
  $T_w$ be the components of $T-\{e_1,e_2\},$ containing $u,v$
  and $w,$ respectively.
  From eigenequations of $M$, we have
  \be \label{u}
 & \lambda x_u = \sum_{y \in T_1\cap V_x} d_{uy}x_y + \sum_{y \in T_2\cap V_x} (d_{vy}+1)x_y + \sum_{y \in T_3\cap V_x} g(d_{wy} +
2)x_y \\
  \label{v}
 & \lambda x_v = \sum_{y \in T_1\cap V_x} (d_{uy}+1)x_y + \sum_{y \in T_2\cap V_x} d_{vy}x_y + \sum_{y \in T_3\cap V_x} (d_{wy} + 1)x_y \\
   \label{w}
 & \lambda x_w = \sum_{y \in T_1\cap V_x} (d_{uy}+2)x_y + \sum_{y \in T_2\cap V_x} (d_{vy}+1)x_y + \sum_{y \in T_3\cap V_x} d_{wy}x_y. \ee

 It follows from \cref{u,v,w} that
\beq
\lambda (x_u+x_w-2x_v)&=\sum_{y \in T_2\cap V_x} 2x_y \\
&\geq 2x_v \\
 &>0 \qquad \text{as }  x_v > 0 \\
\text{ Thus } \hs x_u+x_w& >2x_v \quad \text{as } \lambda > 0.
\eeq

As $T'$ is an arbitrary subtree of the subgraph of $T$ generated by $V_x,$ hence the result follows.
\qed

\bt \label{d}
For a connected graph of diameter $d$, the integers $0,1,\ldots ,d$ are always its distance Pareto eigenvalues.
\et
\pf If $G$ is a graph with diameter $d$, then $\D(P_d)$ is a principal sub-matrix of $\D(G)$ and therefore by \cref{main}  every $\lambda \in \Pi(P_d)$ implies  $\lambda \in\Pi(G).$ Let
\beq A_k=\left(\begin{array}{rr}
0 & k \\
k & 0
\end{array}\right).\eeq
\par
Then for $1\leq k\leq d,$ $A_k$ is a principal sub-matrix of $\D(P_d)$ and $\rho(A_k)=k.$ Therefore $k\in\Pi(P_d)$ for $1\leq k\leq d.$ Besides $0$ being a diagonal element of $\D(P_d)$ is in $\Pi(P_d).$ Hence the result follows. \qed



\par As for a complete graph the distance matrix and the adjacency matrix coincide, therefore from \cite{see18} we have the following.

\bl \label{k} \cite{see18} For any positive integer $n,$ $\Pi(K_n)=\{0, 1, \ldots, n-1\}.$
\el

\bd
If $A$ and $B$ are two nonnegative matrices then we say that $A$ dominates $B$ if either of the following two cases hold
\ben
\item $A$ and $B$ are of same size and upto permutation similarity $A\geq B,$ $A\neq B.$
\item $A$ is permutation similar to $\left(\begin{array}{rr}
B & C \\
D & E
\end{array}\right)$ and at least one of C,D and E is a nonzero matrix.
\een
\ed

From \cref{mi}, \cref{cauchy} and \cref{non} we get the following.
\bl \label{dom}
If $A$ and $B$ are two symmetric nonnegative irreducible matrices, then \\
$A$ dominates $B$ implies $\rho(A)>\rho(B).$
\el

\bl \cite{ah14} \label{srSn} Distance spectral radius of the star $S_n$ is
\beq n-2+\sqrt{(n-2)^2+n-1} \eeq
\el

\bt \label{Sn} There are exactly $2(n-1)$ distance Pareto eigenvalues of $S_n$ and they are
\beq \mu_{2k}=2(k-1), \, \mu_{2k-1}=k-1+\sqrt{k^2-3k+3}
\text{  where }  k=1,\ldots, n-1.
\eeq
\et
\pf Upto permutation similarity there are exactly two distinct principal sub-matrices of $\D(S_n)$ of order $k=2,\ldots,n-1$ and they are $\D(S_k)$ and $2(\J_k-I_k).$ Clearly the later always dominates the former one. Besides $\r(2(\J_k-I_k))=2(k-1)$ and from \cref{srSn}, we have $\r(\D(S_k))=k-2+\sqrt{(k-2)^2+k-1}.$ By routine calculation it can be easily shown that
\be \label{s} \r(\D(S_{k+1}))> \r(2(\J_k-I_k))>\r(\D(S_k)) \hs \text{ for }k=2,\ldots,n-1\ee

As $0$ is always a distance Pareto eigenvalue of a connected graph and $\r(\D(S_n))$ is the largest distance Pareto eigenvalue of $S_n,$ using \cref{s} we get all the distance Pareto eigenvalues of $S_n$ as needed.

\qed

\bt  \label{n}  If $G$ is a connected graph of order $n$ and diameter $d$ then $|\Pi(G)| \geq n+d-1,$ with equality if and only if $G=P_3$ or $K_n.$
\et 
\pf From \cref{d}, we have  $0,1,\ldots, d$ as distance Pareto eigenvalues of $G.$  Let $ A_2=\left(\begin{array}{rr}
0 & d \\
d & 0
\end{array}\right).$  For $n\le 2$ there is nothing to prove. If $n\ge 3$ then for $i=3,\ldots,n$ we can have sub-matrices $A_i$ of $\D(G)$ of order $i$ such that $A_i$ dominates $A_{i-1}.$ Therefore by \cref{dom}, $\r(A_i)>\r(A_{i-1}).$ Thus we get $(d+1)+(n-2)=n+d-1$ distance Pareto eigenvalues of $G$ as follows
\beq 0<1<\cdots <d<\r(A_3)<\cdots <\r(A_n).
\eeq

Hence 
\be \label{nbound} |\Pi(G)| &\geq n+d-1.   \ee 

Now from \cref{k}, we have $|\Pi(K_n)|= n$ and by direct calculation $|\Pi(P_3)|= 4.$ Therefore for $G=P_3,\, K_n$ the equality holds in \cref{nbound}.
\par If $G \ne P_3,\, K_n$ then $n\geq 4$ and $d\ge 2$ and thus $\D(P_3)$ is a principal sub-matrix of $\D(G).$
For $d\ge 3$ we see that $\r(\D(P_3))=1+\sqrt{3}<\r(A_3) $ as minimum row sum of $A_3$ is at least $d+1\ge 4.$ Again if $d=2$ then either $P=\bpm 0 &2 &2 \\ 2 &0 &1 \\ 2 &1 &0  \epm $ or $Q=2(\J_3-I_3)$ must be a principal submatrix of $\D(G).$ So we can choose $A_3=P$ or $A_3=Q,$ whichever be the case. Also $1+\sqrt{3}<\min \{ \r(P),\, \r(Q) \}. $
\par In either situation we get $1+\sqrt{3} \in \Pi (G)$ in addition to the above listed $n+d-1$ distance Pareto eigenvalues of $G.$ Thus $|\Pi(G)| \geq n+d.$

\qed

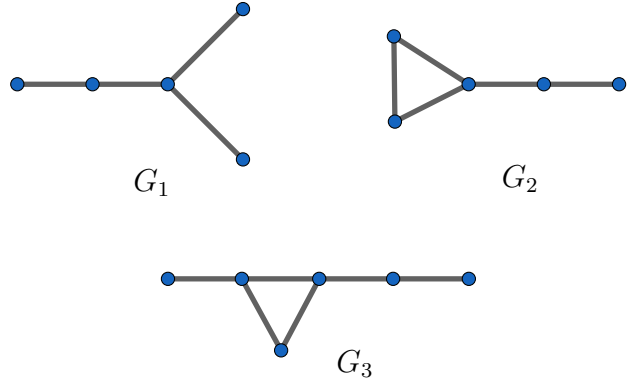
\begin{figure}
    \centering
    \begin{tikzpicture}[line cap=round,line join=round,>=triangle 45,x=1cm,y=1cm]
\clip(-1.7586206896551726,-0.1877483443708616) rectangle (9,6);
\draw [line width=2pt,color=wrwrwr] (2.0027685020151527,1.4110914217620147)-- (6.002768502015153,1.4110914217620147);
\draw [line width=2pt,color=wrwrwr] (0,4)-- (1,4);
\draw [line width=2pt,color=wrwrwr] (1,4)-- (2,4);
\draw [line width=2pt,color=wrwrwr] (2,4)-- (3,5);
\draw [line width=2pt,color=wrwrwr] (2,4)-- (3,3);
\draw [line width=2pt,color=wrwrwr] (5.006249573116593,4.636961956150539)-- (5.019909842223894,3.5031596202445217);
\draw [line width=2pt,color=wrwrwr] (5.019909842223894,3.5031596202445217)-- (6,4);
\draw [line width=2pt,color=wrwrwr] (6,4)-- (5.006249573116593,4.636961956150539);
\draw [line width=2pt,color=wrwrwr] (6,4)-- (7,4);
\draw [line width=2pt,color=wrwrwr] (7,4)-- (8,4);
\draw [line width=2pt,color=wrwrwr] (3.507027411067685,0.4594174180757188)-- (2.985141667110684,1.4110914217620147);
\draw [line width=2pt,color=wrwrwr] (3.507027411067685,0.4594174180757188)-- (4.013563574320068,1.4110914217620147);
\draw (1.3985801217038543,3.013907284768212) node[anchor=north west] {$G_1$};
\draw (6.296146044624748,3.068543046357616) node[anchor=north west] {$G_2$};
\draw (4.098377281947262,0.5990066225165558) node[anchor=north west] {$G_3$};
\begin{scriptsize}
\draw [fill=rvwvcq] (0,4) circle (2.5pt);
\draw [fill=rvwvcq] (1,4) circle (2.5pt);
\draw [fill=rvwvcq] (2,4) circle (2.5pt);
\draw [fill=rvwvcq] (3,5) circle (2.5pt);
\draw [fill=rvwvcq] (3,3) circle (2.5pt);
\draw [fill=rvwvcq] (5.006249573116593,4.636961956150539) circle (2.5pt);
\draw [fill=rvwvcq] (5.019909842223894,3.5031596202445217) circle (2.5pt);
\draw [fill=rvwvcq] (6,4) circle (2.5pt);
\draw [fill=rvwvcq] (7,4) circle (2.5pt);
\draw [fill=rvwvcq] (8,4) circle (2.5pt);
\draw [fill=rvwvcq] (2.0027685020151527,1.4110914217620147) circle (2.5pt);
\draw [fill=rvwvcq] (6.002768502015153,1.4110914217620147) circle (2.5pt);
\draw [fill=rvwvcq] (2.985141667110684,1.4110914217620147) circle (2.5pt);
\draw [fill=rvwvcq] (4.013563574320068,1.4110914217620147) circle (2.5pt);
\draw [fill=rvwvcq] (4.995936739415599,1.4110914217620147) circle (2.5pt);
\draw [fill=rvwvcq] (3.507027411067685,0.4594174180757188) circle (2.5pt);
\end{scriptsize}
\end{tikzpicture}

\caption{Graphs with maximum number of distance Pareto eigenvalues of order 5 and 6}
    \label{fig3}
\end{figure}

{\bf Note:} From \cref{n} we see that for any connected graph $G$ with $n$ vertices, $|\Pi(G)|\geq n$
and the equality is achieved if and only if $G=K_n.$
But the exact number of distance Pareto eigenvalues of a connected graph of order $n$ is not known. Now for any $n \times n$ matrix there are exactly $2^n-1$ principal sub-matrices, so by \cref{main} we can not have more than $2^n-1$ distance Pareto eigenvalues of any connected graph with $n$ vertices. Also we can see that for any connected graph with $n\ge 2$ vertices, all the principal sub-matrices of $\D(G)$ are not distinct. Besides we may have some principal sub-matrices which are distinct but having same spectral radius. Therefore for any connected graph $G$ of order $n$ if there are $s(G)$ distinct principal sub-matrices of $\D(G)$ then $|\Pi(G)|\le s(G),$ with equality if and only if all the distinct sub-matrices of $\D(G)$ have distinct spectral radii. In this regard another question arises that among all connected graphs of given order which graph(s) will have maximum number of distance Pareto eigenvalues. By direct calculation we have seen that among all connected graphs of order $n=2,3,4$ the path graph $P_n$ has the maximum number of distance Pareto eigenvalues. The numbers are respectively 2, 4 and 7. But for graphs of order 5, we see that the path graph $P_5$ together with graphs $G_1$ and $G_2$ of \cref{fig3} have maximum (here 13) number of distance Pareto eigenvalues. Again for the class of graphs of order 6, the graph $G_3$ in \cref{fig3} attains uniquely the maximum (30) number of distance Pareto eigenvalues. Therefore this seems to be an interesting problem to characterize all graphs  for which the maximum number of distance Pareto eigenvalues occur among all connected graphs of given order.  We here pose a question which perhaps require deep investigation. Can distance Pareto eigenvalues of a connected graph uniquely determine the graph? In addition to that it is worth studying how fast distance Pareto eigenvalues grow when number of vertices increases.
We leave these problems for future research scope.

\bt \label{bound} If $G$ is a graph with $n$ vertices, then
\beq \r_k(G)\geq n-k \hs \text{for} \hs k=1,2,\ldots,n. \eeq
Equality holds if and only if $G=K_n.$
\et
\pf If $G=K_n,$ then using \cref{k} we are done.
\par If $G\neq K_n,$ then for $i=1,2,\ldots,n-2,$ we can have $A_i$ of order $n-i+1$ as principal sub-matrix of $\D(G)$ such that $A_{1}=\D(G),$ $A_{n-2}=\D(P_3),$ $A_{i+1}$ is a principal sub-matrix of  $A_i$ for $i=1,2,\ldots, n-3.$

\par Besides $A_i$ dominates $\J_{n-i+1}-I_{n-i+1}$ for $i=1,2,\ldots, n-3.$ \\
But \beq\r(\J_{n-i+1}-I_{n-i+1})=n-i.\eeq
\be \label{n-k}
\text{ Therefore }\quad \r_i(G) > n-i \hs \text{for} \hs i=1,2,\ldots, n-3.
\ee
Now since $A_{n-2}=\D(P_3)$  and  $\Pi(P_3)=\{0,1,2,1+\sqrt{3}\},$ we get
\be \label{01}
 \r_{n-2}(G)\geq 1+\sqrt{3}>2, \, \r_{n-1}(G)\geq 2>1 \, \text{ and }\, \r_n (G)\geq 1>0.
\ee

Combining \cref{n-k,01}, we get the result as desired.
\qed

\section{Second largest distance Pareto eigenvalue}\label{sec5}

The largest distance Pareto eigenvalue of a connected graph is nothing but the distance spectral radius of the graph. Also in last few decades distance spectral radius have been extensively studied. In this section we study some bounds for the second largest distance Pareto eigenvalue.

\bt \label{r2def} If $G$ is a connected graph with at least two vertices, then
\beq \r_2(G)=\max \{ \r(A): A\in P\},
\eeq

where $P=\{ (\D(G))(v): v\in V(G), d_v>1 \}$
\et
\pf Since $\D(G)$ dominates $A$ for every $A\in P,$ we have
\beq \r_2(G)\leq\r(A)<\r_1(G) \text{  for every } A\in P. \eeq

Also for every principal sub-matrix $B$ of $\D(G)$ of order $n-2$ or less, there exist $A\in P$ which dominates $B$ and therefore $\r(B)<\r(A).$ Thus $\r_2(G)$ must be equal to the largest $\r(A)$ for all possible $A\in P.$ \\

Now if $u$ is a pendent vertex in $G$ and $v$ is a quasi-pendent vertex in $G$ adjacent to $u,$ then the principal sub-matrix of $\D(G)$ obtained by removing row and column corresponding to vertex $v$ dominates the principal sub-matrix of $\D(G)$ obtained by removing row and column corresponding to vertex $u.$ Therefore in calculating $\r_2(G)$ we can ignore those  principal sub-matrices of $\D(G)$ which are obtained by removing row and column corresponding to pendent vertex. Hence the result follows.\qed

\bt If $G$ is a graph of order $n$ with a vertex of degree $n-1,$ then
\beq n-2\leq \r_2(G)\leq 2(n-2), \eeq
the right hand equality holds if and only if $G=S_n$ and the left hand equality holds if and only if $G=K_n.$
\et
\pf From \cref{bound}, we see that $\r_2(G)\geq n-2$ and equality holds if and only if $G=K_n.$
\par Now suppose $G\neq K_n.$ Let $v\in V(G)$ with $d_G(v)=n-1.$ and $A$ be the principal sub-matrix of $\D(G)$ obtained by deleting row and column of $\D(G)$ corresponding to vertex $v.$ Then Clearly $\r_2(G)=\r(A),$ as $A$ dominates any other principal sub-matrix of $\D(G)$ of order $n-1.$
\par From \cref{Sn}, $\r_2(S_n)=2(n-1).$ Now if $G\neq S_n,$ then we can find two vertices $u,w(\neq v)\in V(G)$ such that $u \sim w$ in $G.$ Thus $2(\J_{n-1}-I_{n-1})$ dominates $A.$ Therefore $\r_2(G)=\r(A)<2(n-1).$
\qed

\bt \label{thmsn} If $G$ is a connected graph of order $n$ and diameter $2$ then $\r_2(G)\le 2(n-2),$ with equality if and only if $G=S_n.$
\et 
\pf By \cref{Sn}, we have $\r_2(S_n)=2(n-2).$ Now if $G\ne S_n,$ then there is at least two non pendent vertices $u,v \in V(G).$ So there are at least two $1's$ in each of $u-th$ and $v-th$ rows(columns) of $\D(G).$ Therefore if $A$ is any principal submatrix of $\D(G)$ of order $n-1,$ then $2(\J_{n-1}-I_{n-1})$ dominates $A.$ Thus $2(n-2)>\r_2(G).$
\qed 

\bt If $G$ is a connected graph with $n$ vertices and $\omega(G)\geq n-1,$ then
\beq n-2\leq\r_2(G)\leq\frac{n-3+\sqrt{n^2+10n-23}}{2} ,\eeq
with equality in the left hand side if and only if $G=K_n$ and equality in the right hand side if and only if $G=K_{n-1}^1.$
\et

\pf From \cref{k} we have $\r_2(G)\geq n-2$ with equality if and only if $G=K_n.$\\
Now as $\omega(G)\geq n-1,$ we can take $H=K_{n-1}$ to be a subgraph of $G$ and $v\in V(G)- V(H).$ Then upto permutation similarity we have
\beq \D(G)=\bpm
\D(K_{n-1}) & x \\
x^T & 0
\epm, \hs
\text{ where } x_u &=1  \text{ if } u\sim v\\
&=2 \text{ otherwise.}
\eeq
Since $G$ is connected, there exists $w\in V(H)$ with $w\sim v.$ If $B(G)$ is the principal sub-matrix of $\D(G)$  obtained by deleting row and column corresponding to $w$ then clearly $\r_2(G)=\r(B).$ and it can be easily observed that $B(K_{n-1}^p)$ dominates $B(K_{n-1}^{p-1})$  for $p=2,\ldots, n.$ Therefore
\beq \r_2(K_{n-1}^1)> \r_2(K_{n-1}^2)>\cdots >\r_2(K_{n-1}^{n-1})>\r_2(K_{n-1}^n)
\eeq

Thus we get $\r_2(G)\leq \r(\D(K_{n-1}^1)),$ with equality if and only if $G=K_{n-1}^1.$\\
\par Now let $(\r,y)$ be the eigenpair of $B(K_{n-1}^1).$ Then due to symmetry $y_2=y_3=\cdots =y_{n-1}.$\\
\par From eigenequations we have
\be\label{y1} 2(n-1)y_2 &=\r y_1\\
 \label{y2} 2y_1+(n-3)y_2&=\r y_2 \ee

From \cref{y1,y2} we get $\r^2-(n-3)\r-4(n-2)=0.$ Which gives
\beq \r=\frac{n-3+\sqrt{n^2+10n-23}}{2}. \eeq
\par Hence the result follows.
\qed

\bt \label{kn-e} For any non complete connected graph $G$ with $n$ vertices,
\beq \r_2(G)\geq \frac{n-2+\sqrt{n^2-4n+12}}{2}, \eeq
with equality if and only if $G=K_n-e.$
\et
\pf First we find the expression for $\r(\D(K_n-e).$ If $(\r, x)$ is the eigenpair for $\D(K_n-e),$ where $e=\{1,2\}.$ Then due to symmetry we have
\be \label{341}x_1=x_2 \text{ and } x_3=\cdots =x_n. \ee

Also we have

\be \label{342}  \D(K_n-e)=\bpm
0 & 2 & \1_{n-2}^T \\
2 & 0 & 1 \\
\1_{n-2} & 1 & \D(K_{n-2})
\epm \ee
 Using \cref{341,342} in eigenequations of $\D(K_n-e)$, we get

\be \label{343} 2x_1+(n-2)x_3&=\r x_1 \\
\label{344} 2x_1+(n-3)x_3&=\r x_3
\ee

From \cref{343,344} we have $\r^2-(n-1)\r-2=0.$ Thus we get
\be \label{345} \r(\D(K_n-e))=\frac{n-1+\sqrt{(n-1)^2+8}}{2}. \ee

Clearly $\D(K_{n-1}-e)$ and $\J_{n-1}-I_{n-1}$ are the only two distinct sub-matrices of $\D(K_n-e)$ of order $n-1$ and the former dominates the later one. Therefore
\beq  \hs \r_2(K_n-e))&=\r (\D(K_{n-1}-e)) \\
&=\frac{n-2+\sqrt{n^2-4n+12}}{2} \quad \text{ [using \cref{345}]}\eeq
\par Now if $G\neq K_n-e,$ then any principal sub-matrix of $\D(G)$ of order $n-1$ dominates $\D(K_{n-1}-e)$  and therefore $\r_2(G)> \frac{n-2+\sqrt{n^2-4n+12}}{2}.$
\par Hence the result follows.
\qed

\bc For any non complete connected graph $G$ with n vertices,
\[ \r_2(G)\geq n-2+\frac{2}{n-1} \]
\par equality holds if and only if $G=P_3.$
\ec
\pf As $G$ is connected and non complete, therefore $n \geq 3.$ Which implies that
\beq
n^2-4n+12 & \geq \Big(n-2+\frac{4}{n-1}\Big)^2 .\eeq

Thus
\be \label{p3} \frac {n-2+\sqrt{n^2-4n+12}}{2} & \geq n-2+\frac{2}{n-1} \ee
It can easily be shown that equality in \cref{p3} holds if and only if $n=3.$ Therefore by \cref{kn-e} we get
\[ \r_2(G)\geq n-2+\frac{2}{n-1} ,\]
\par with equality if and only if $G=K_3-e=P_3.$
\qed


\bt \label{W-T} If $G$ is a connected graph of order $n$ so that minimum transmission occur at a vertex $v\in G$ and $\x$ is the normalized distance Pareto eigenvector corresponding to $\r_2,$ then
\[ \r_2(G)\geq \frac{2[W- Tr(v)]}{n-1},\]

with equality if and only if $x_u=\frac{1}{\sqrt{n-1}}$ for $u\ne v. $
\et
\pf Let $D_v$ be the sub-matrix of $\D(G)$ obtained by deleting row and column corresponding to vertex $v.$ Then the average row sum of $D_v$ equals $\frac {2(W-Tr(v))}{n-1}.$ Therefore using \cref{avg} we get
\beq \r_2(G)&\ge \r(D_v) \\ &\geq \frac {2(W-Tr(v))}{n-1}.
\eeq

Now equality holds in the above expression if and only if $\r_2(G)= \r(D_v)$ and all the row sums of $D_v$ are equal which is the case if and only if $x_v=0$ and $x_u$ is constant for $u\ne v.$ Thus we get the required condition from the facts that $\x$ is normalized and exactly one component of $\x$ is zero.
\qed \\
{\bf Note:} The equality in \cref{W-T} holds for several graphs like complete graph, star graph, wheel graph etc.

\bt For any connected graph $G$ of order $n$ other than $K_n$ and $K_n-e$
\beq \r_2(G)\geq \frac{n-2+\sqrt{n^2-4n+20}}{2}, \eeq
with equality if and only if $G=K_n-\{e_1, e_2\},$ where $e_1$ and $e_2$ are not incident in $K_n.$
\et
\pf Let $G_1=K_n-\{e_1, e_2\},$ where $e_1,\, e_2\in E(K_n)$ are not incident and $G_2=K_n-\{f_1, f_2\},$ where $f_1, \, f_2 \in E(K_n)$ are incident.
Upto permutation similarity there are exactly two distinct principal sub-matrices of $\D(G_1)$ of order $n-1$ and they are given by

\beq
&M=\bpm
2(\J_2-I_2) & \J_2-I_2 & \J-I \\
\J_2-I_2 & 2(\J_2-I_2) & \J_2-I_2 \\
\J-I & \J_2-I_2 & \J-I
\epm \\
\\
\text{and } \hs &N=\bpm
2(\J_2-I_2) & \J-I \\
\J-I & \J-I
\epm .
\eeq

Clearly $M$ dominates $N$ and therefore $\r_2(G_1)=\r(M).$\\
If $(\r,x)$ be the eigenpair of $M$, then due to symmetry

\be \label{351}  x_1=x_2=x_3=x_4=a \text{ (say) } \\
\label{352} \text{and } \quad x_5=\cdots =x_n=b \text{ (say) }\ee

Using \cref{351,352} in eigenequations of $M$ we get

\be \label{353} 2a+2a+(n-5)b=\r a \\
\label{354} 4a+(n-6)b=\r b
\ee

From \cref{353,354} we get\\
\[ \r_2(G_1)=\r(M)=\frac{n-2+\sqrt{(n-2)^2+16}}{2}. \]

Now if $G\neq G_1,$ then we consider the following two cases

\bdsc
\item{Case I:} $G$ has at most $\binom{n}{2}-3$ edges. In this case there are principal sub-matrices of $\D(G)$ of order $n-1$ which dominates $M.$ Therefore we get

\[ \r_2(G)>\r (M)=\frac{n-2+\sqrt{(n-2)^2+16}}{2}.\]

\item{Case II:} If $G=G_2,$ then it can be observed that $\r_2(G_2)=\r (A),$ where
\[ A= \begin{pmatrix}
&0 &2 &2 & \1\\
&2 &0 &1 & 1 \\
&2 &1 &0 & \1 \\
& \1 &\1 & \1 & \J-I
\end{pmatrix}.
\]

If $(\r , x)$ is the eigenpair of $A,$ then due to symmetry,
\[ x_2=x_3 \text{ and } x_4=x_5=\cdots =x_n. \]

Therefore from eigenequations we have

\be
\label{3511} 2x_2+2x_2+(n-4)x_4=\r x_1 \\
\label{3512} 2x_1+x_2+(n-4)x_4=\r x_2 \\
\label{3513} x_1+2x_2+(n-5)x_4=\r x_4
\ee
From \crefrange{3511}{3513}, we have

\beq \begin{vmatrix}
\r &-4 &-(n-4) \\
-2 &\r-1 &-(n-4)\\
-1 & -2 & \r-n+5
\end{vmatrix} =0\\
\eeq
Which implies $(\r+2)[\r^2-(n-1)\r +(n-8)]=0$

Therefore $\r (A)$ is the largest root of $y^2-(n-1)y+n-8=0.$

\par Let \[ f(y)= y^2-(n-1)y+n-8.\]
Then it can be verified that
\[ f\bigg( \frac{n-2+\sqrt{(n-2)^2+16}}{2}\bigg)=\frac{n-\sqrt{n^2-4n+20}}{2}-3.\]
Now $G\neq K_n, K_n-e$ gives  $n  \geq 4.$ Which in turn implies that
\beq \frac{n-2+\sqrt{(n-2)^2+16}}{2} & < 3 \\
\text{ Thus } \hs f\bigg( \frac{n-2+\sqrt{(n-2)^2+16}}{2}\bigg) & < 0 \quad \text{ for all } \,  n \geq 3.
\eeq

Hence the largest root of $f(y)=0$  must be greater than $\frac{n-2+\sqrt{(n-2)^2+16}}{2}.$
\edsc
i.e. \[ \r_2(G) \geq \frac{n-2+\sqrt{(n-2)^2+16}}{2}, \] with equality if and only if $G=G_1.$
\qed

\bt \label{2>2} If $\lambda_2(G)$ is the second largest distance eigenvalue of a connected graph $G,$ then $\r_2(G)>\lambda_2(G).$
\et
\pf Let \[\ D(G)=\bpm
0 & \x^T \\
\x & B
\epm  \text{ such that } \r_2(G)=\r(B). \]

From Cauchy's inequality (\cref{cauchy}), we have \\
\[ \r_2(G)\geq \lambda_2(G),\]
\par with equality if and only if there exists $z\in \mathbb{C}^{n-1}-\{ \textbf{0} \}$ such that
\[ Bz=\r(B)z \text{ and } \x^Tz=0.\]
 Here $\x$ is always a positive vector and $z$ being an eigenvector corresponding to the perron value of a nonnegative irreducible matrix is real and is either positive or negative. In either case $\x^Tz=0$ can never hold. Hence $\r_2(G)>\lambda_2(G).$
\qed

\bt \label{G-e} If $G$ and $G'=G-e$ are connected graphs then $\r_2(G')\geq \r_2(G).$
\et
\pf Suppose $\r_2(G)=\r(A),$ where $A$ is the sub-matrix of $\D(G)$ obtained by deleting row and column of $\D(G)$ corresponding to vertex $v\in G.$ Let $B$ be the sub-matrix of $\D(G')$ obtained by deleting row and column corresponding to vertex $v.$

 Then clearly either $B=A$ or $B$ dominates $A.$ Thus $\r_2(G')\geq \r_2(G).$ \qed \\

\begin{figure}
\centering

\begin{tikzpicture}
[line cap=round,line join=round,>=triangle 45,x=1cm,y=1cm]
\clip(-9.813740993914806,-2.0456953642384113) rectangle (-3.924970993914806,2.8543046357615887);
\fill[line width=2pt,color=rvwvcq,fill=rvwvcq,fill opacity=0.10000000149011612] (-8.26,0.9) -- (-8.86,0.42) -- (-8.588903068803292,-0.29896206707706696) -- (-7.821355951077937,-0.2633050611525747) -- (-7.618082675553357,0.4776942475228838) -- cycle;
\fill[line width=2pt,color=rvwvcq,fill=rvwvcq,fill opacity=0.10000000149011612] (-5.6148167551344885,1.0019671538621158) -- (-6.214816755134486,0.5219671538621156) -- (-5.94371982393778,-0.19699491321494922) -- (-5.1761727062124265,-0.16133790729045694) -- (-4.972899430687847,0.5796614013850003) -- cycle;
\draw [line width=2pt,color=rvwvcq] (-8.26,0.9)-- (-8.86,0.42);
\draw [line width=2pt,color=rvwvcq] (-8.86,0.42)-- (-8.588903068803292,-0.29896206707706696);
\draw [line width=2pt,color=rvwvcq] (-8.588903068803292,-0.29896206707706696)-- (-7.821355951077937,-0.2633050611525747);
\draw [line width=2pt,color=rvwvcq] (-7.821355951077937,-0.2633050611525747)-- (-7.618082675553357,0.4776942475228838);
\draw [line width=2pt,color=rvwvcq] (-7.618082675553357,0.4776942475228838)-- (-8.26,0.9);
\draw [line width=2pt,color=wrwrwr] (-8.238866441092208,1.7748557772081432)-- (-8.86,0.42);
\draw [line width=2pt,color=wrwrwr] (-8.238866441092208,1.7748557772081432)-- (-8.26,0.9);
\draw [line width=2pt,color=wrwrwr] (-8.238866441092208,1.7748557772081432)-- (-8.588903068803292,-0.29896206707706696);
\draw [line width=2pt,color=wrwrwr] (-8.238866441092208,1.7748557772081432)-- (-7.821355951077937,-0.2633050611525747);
\draw [line width=2pt,color=wrwrwr] (-8.238866441092208,1.7748557772081432)-- (-7.618082675553357,0.4776942475228838);
\draw [line width=2pt,color=rvwvcq] (-5.6148167551344885,1.0019671538621158)-- (-6.214816755134486,0.5219671538621156);
\draw [line width=2pt,color=rvwvcq] (-6.214816755134486,0.5219671538621156)-- (-5.94371982393778,-0.19699491321494922);
\draw [line width=2pt,color=rvwvcq] (-5.94371982393778,-0.19699491321494922)-- (-5.1761727062124265,-0.16133790729045694);
\draw [line width=2pt,color=rvwvcq] (-5.1761727062124265,-0.16133790729045694)-- (-4.972899430687847,0.5796614013850003);
\draw [line width=2pt,color=rvwvcq] (-4.972899430687847,0.5796614013850003)-- (-5.6148167551344885,1.0019671538621158);
\draw [line width=2pt,color=wrwrwr] (-5.593683196226697,1.8768229310702595)-- (-5.6148167551344885,1.0019671538621158);
\draw [line width=2pt,color=wrwrwr] (-5.593683196226697,1.8768229310702595)-- (-5.94371982393778,-0.19699491321494922);
\draw [line width=2pt,color=wrwrwr] (-5.593683196226697,1.8768229310702595)-- (-5.1761727062124265,-0.16133790729045694);
\draw [line width=2pt,color=wrwrwr] (-5.593683196226697,1.8768229310702595)-- (-4.972899430687847,0.5796614013850003);
\draw (-8.368424219066934,-0.4329470198675503) node[anchor=north west] {$G_1$};
\draw (-5.794326987829613,-0.30546357615894104) node[anchor=north west] {$G_2$};
\draw (-9.134270121703852,-0.9698675496688748) node[anchor=north west] {$\rho _2(G_1)=6$};
\draw (-6.25420050709939,-0.9061258278145701) node[anchor=north west] {$\rho_2(G_2)=6$};

\begin{scriptsize}
\draw [fill=rvwvcq] (-8.26,0.9) circle (2.5pt);
\draw [fill=rvwvcq] (-8.86,0.42) circle (2.5pt);
\draw [fill=rvwvcq] (-8.588903068803292,-0.29896206707706696) circle (2.5pt);
\draw [fill=rvwvcq] (-7.821355951077937,-0.2633050611525747) circle (2.5pt);
\draw [fill=rvwvcq] (-7.618082675553357,0.4776942475228838) circle (2.5pt);
\draw [fill=rvwvcq] (-8.238866441092208,1.7748557772081432) circle (2.5pt);
\draw [fill=rvwvcq] (-5.6148167551344885,1.0019671538621158) circle (2.5pt);
\draw [fill=rvwvcq] (-6.214816755134486,0.5219671538621156) circle (2.5pt);
\draw [fill=rvwvcq] (-5.94371982393778,-0.19699491321494922) circle (2.5pt);
\draw [fill=rvwvcq] (-5.1761727062124265,-0.16133790729045694) circle (2.5pt);
\draw [fill=rvwvcq] (-4.972899430687847,0.5796614013850003) circle (2.5pt);
\draw [fill=rvwvcq] (-5.593683196226697,1.8768229310702595) circle (2.5pt);

\end{scriptsize}

\end{tikzpicture}
\caption{Graphs $G_1$ and $G_2=G_1-e$ with $\r_2(G_1)=\r_2(G_2).$}
\label{fig1}
\end{figure}
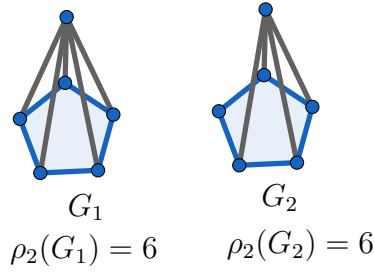

{\bf Note:} The inequality in the \cref{G-e} is not always strict. For example we can consider the graphs $G_1$ and $G_2=G_1-e$ as in \cref{fig1} with $\r_2(G_1)=6=\r_2(G_2).$ In \cref{G-e} if $e$ is not incident with $v$ then clearly $B$ dominates $A$ and therefore $\r_2(G')> \r_2(G).$ Besides if $e$ connects $v$ to $u \in V(G)$ and for some $i, j(\ne v) \in V(G),$ if $d_{ij}$ increases from $G$ to $G'$ then again $B$ dominates $A$ and therefore $\r_2(G')> \r_2(G).$ In this regard the problem of classifying the edges (if any) in a graph $G$  whose removal do not increase $\r_2(G)$ seems to be an interesting problem.

From \cref{G-e} we get the following result as immediate corollary.

\bc \label{r2mtree} Among all connected graphs of given order, second largest distance Pareto eigenvalue is maximum for some tree.
\ec

\bl \label{r2kab} If $a\le b,$ then $\r_2(K_{a,b})=a+b-3+\sqrt{a^2+b^2+b-ab-2a+1}. $
\el
\pf If $a=1$ then by \cref{Sn} we know that $\r_2(K_{a,b})=2(b-1).$  Now if $a\ge 2$ then upto permutation similarity $\D(K_{a,b})$ has exactly two distinct principal submatrices of order $n-1$ namely $\D(K_{a-1,b})$ and $\D(K_{a,b-1}).$ It can be observed that
\beq \r_1(K_{a,b})=a+b-2+\sqrt{a^2+b^2-ab}. \eeq 
Therefore we get
\beq \r(\D(K_{a-1,b}))=a+b-3+\sqrt{a^2+b^2+b-ab-2a+1} \\
\text{and } \quad \r(\D(K_{a,b-1}))=a+b-3+\sqrt{a^2+b^2+a-ab-2b+1}\eeq
Now $a\le b$ implies that $\r(\D(K_{a-1,b}))>\r(\D(K_{a,b-1}))$ and the lemma follows. \qed

\bt If $G$ is a connected bipartite graph of order $n,$ then
\beq \r_2(G)\ge n-3+\sqrt{n^2+n+1+3\lfloor \frac{n}{2}\rfloor(\lfloor\frac{n}{2}\rfloor -n-1 ) } \eeq
equality is attained if and only if $G=K_{\lfloor\frac{n}{2}\rfloor ,\lceil \frac{n}{2}\rceil }$ \et 

\pf Let $G$ be a bipartite graph with $V(G)=V_1\cup V_2$ as vertex bipartition such that $|V_1|=p$ and $|V_2|=q.$ Then from \cref{G-e} we have 

\be \label{edr2i} \r_2(G)\ge \r_2(K_{p,q}). \ee
Now from the proof of \cref{r2kab} we see that for any edge $e=(u,v)$ in $K_{p,q},$ where $p\ge 2$ we can have a vertex $w$ in $K_{p,q}$ different from $u$ and $v$ so that $\r_2(K_{p,q})=\r(A)$ where $A=(\D(K_{p,q}))(w).$ Since $e$ is not incident with $w,$ therefore $(\D(K_{p,q}-e))(w)$ dominates $A$ and therefore $\r_2(K_{p,q}-e))>\r(K_{p,q}).$ Thus the equality in \cref{edr2i} holds if and only if $G=K_{p,q}.$ Again for $p\le \lfloor\frac{n}{2}\rfloor $ writin $q=n-p$ we get from \cref{r2kab}
\beq \r_2(K_{p,n-p})=n-3+\sqrt{n^2+n+1+3p(p-n-1)} \eeq
which is a strictly decreasing function for $p\le \lfloor\frac{n}{2}\rfloor .$ Thus we get
\beq \r_2(K_{1,n-1})>\r_2(K_{2,n-2})>\cdots > \r_2(K_{\lfloor\frac{n}{2}, \lceil\frac{n}{2}\rceil } )   \eeq

Considering all the above arguments we can say that $\r_2(G)\ge \r_2(K_{\lfloor\frac{n}{2}\rfloor ,\lceil \frac{n}{2}\rceil }) $ and equality holds if and only if $G=K_{\lfloor\frac{n}{2}\rfloor ,\lceil \frac{n}{2}\rceil }.$ This proves the theorem.

\qed 


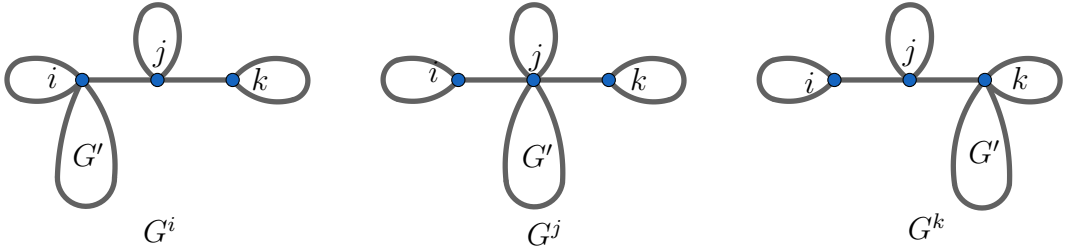
\begin{figure}
    \centering
\begin{tikzpicture}[line cap=round,line join=round,>=triangle 45,x=1cm,y=1cm]
\clip(-0.58378,0.5070860139510065) rectangle (32.23696960176991,4.542022080594171);
\draw [shift={(0.3846438471702293,2.4871806418906877)},line width=2pt,color=wrwrwr]  plot[domain=0.6947599161342964:1.826877887033421,variable=\t]({1*0.8010286442300945*cos(\t r)+0*0.8010286442300945*sin(\t r)},{0*0.8010286442300945*cos(\t r)+1*0.8010286442300945*sin(\t r)});
\draw [shift={(0.39579926378568586,3.5011563388713425)},line width=2pt,color=wrwrwr]  plot[domain=4.452523221142281:5.590741217090987,variable=\t]({1*0.7849943984722739*cos(\t r)+0*0.7849943984722739*sin(\t r)},{0*0.7849943984722739*cos(\t r)+1*0.7849943984722739*sin(\t r)});
\draw [shift={(0.2735871948022161,3.0043384281235577)},line width=2pt,color=wrwrwr]  plot[domain=1.9130770913733177:4.417617444885296,variable=\t]({1*0.27362159110408113*cos(\t r)+0*0.27362159110408113*sin(\t r)},{0*0.27362159110408113*cos(\t r)+1*0.27362159110408113*sin(\t r)});
\draw [shift={(2.3722366472923073,3.5456369524222224)},line width=2pt,color=wrwrwr]  plot[domain=2.7427875861151514:4.113705584423451,variable=\t]({1*0.6605148033435952*cos(\t r)+0*0.6605148033435952*sin(\t r)},{0*0.6605148033435952*cos(\t r)+1*0.6605148033435952*sin(\t r)});
\draw [shift={(1.575071066600976,3.576017594944613)},line width=2pt,color=wrwrwr]  plot[domain=-0.9352129575633459:0.3148685697333037,variable=\t]({1*0.7157938726514836*cos(\t r)+0*0.7157938726514836*sin(\t r)},{0*0.7157938726514836*cos(\t r)+1*0.7157938726514836*sin(\t r)});
\draw [shift={(2.0091549289237203,3.748856535752751)},line width=2pt,color=wrwrwr]  plot[domain=0.19557120262823133:2.9280039203807275,variable=\t]({1*0.2513102710947318*cos(\t r)+0*0.2513102710947318*sin(\t r)},{0*0.2513102710947318*cos(\t r)+1*0.2513102710947318*sin(\t r)});
\draw [shift={(3.571244562288569,2.571989219789647)},line width=2pt,color=wrwrwr]  plot[domain=1.2413327694755845:2.4985652239452047,variable=\t]({1*0.7138021980356563*cos(\t r)+0*0.7138021980356563*sin(\t r)},{0*0.7138021980356563*cos(\t r)+1*0.7138021980356563*sin(\t r)});
\draw [shift={(3.597867360796995,3.467597203185726)},line width=2pt,color=wrwrwr]  plot[domain=3.8053285983434444:4.992084661918382,variable=\t]({1*0.7590075925400728*cos(\t r)+0*0.7590075925400728*sin(\t r)},{0*0.7590075925400728*cos(\t r)+1*0.7590075925400728*sin(\t r)});
\draw [shift={(3.7364329193972488,2.9920425106768636)},line width=2pt,color=wrwrwr]  plot[domain=-1.2982974755783836:1.3187811236619362,variable=\t]({1*0.2636871775680132*cos(\t r)+0*0.2636871775680132*sin(\t r)},{0*0.2636871775680132*cos(\t r)+1*0.2636871775680132*sin(\t r)});
\draw [line width=2pt,color=wrwrwr] (1,3)-- (2,3);
\draw [line width=2pt,color=wrwrwr] (2,3)-- (3,3);
\draw [shift={(5.384643847170229,2.4871806418906846)},line width=2pt,color=wrwrwr]  plot[domain=0.694759916134299:1.826877887033419,variable=\t]({1*0.8010286442300969*cos(\t r)+0*0.8010286442300969*sin(\t r)},{0*0.8010286442300969*cos(\t r)+1*0.8010286442300969*sin(\t r)});
\draw [shift={(5.395799263785686,3.501156338871341)},line width=2pt,color=wrwrwr]  plot[domain=4.45252322114228:5.590741217090989,variable=\t]({1*0.7849943984722729*cos(\t r)+0*0.7849943984722729*sin(\t r)},{0*0.7849943984722729*cos(\t r)+1*0.7849943984722729*sin(\t r)});
\draw [shift={(7.372236647292307,3.545636952422223)},line width=2pt,color=wrwrwr]  plot[domain=2.742787586115152:4.113705584423452,variable=\t]({1*0.6605148033435952*cos(\t r)+0*0.6605148033435952*sin(\t r)},{0*0.6605148033435952*cos(\t r)+1*0.6605148033435952*sin(\t r)});
\draw [shift={(6.575071066600976,3.576017594944613)},line width=2pt,color=wrwrwr]  plot[domain=-0.9352129575633459:0.31486856973330346,variable=\t]({1*0.7157938726514839*cos(\t r)+0*0.7157938726514839*sin(\t r)},{0*0.7157938726514839*cos(\t r)+1*0.7157938726514839*sin(\t r)});
\draw [shift={(7.009154928923721,3.7488565357527506)},line width=2pt,color=wrwrwr]  plot[domain=0.19557120262823308:2.928003920380726,variable=\t]({1*0.25131027109473236*cos(\t r)+0*0.25131027109473236*sin(\t r)},{0*0.25131027109473236*cos(\t r)+1*0.25131027109473236*sin(\t r)});
\draw [shift={(8.571244562288566,2.5719892197896512)},line width=2pt,color=wrwrwr]  plot[domain=1.2413327694755818:2.498565223945207,variable=\t]({1*0.713802198035651*cos(\t r)+0*0.713802198035651*sin(\t r)},{0*0.713802198035651*cos(\t r)+1*0.713802198035651*sin(\t r)});
\draw [shift={(8.597867360796995,3.4675972031857247)},line width=2pt,color=wrwrwr]  plot[domain=3.805328598343443:4.99208466191838,variable=\t]({1*0.759007592540072*cos(\t r)+0*0.759007592540072*sin(\t r)},{0*0.759007592540072*cos(\t r)+1*0.759007592540072*sin(\t r)});
\draw [shift={(8.736432919397249,2.992042510676863)},line width=2pt,color=wrwrwr]  plot[domain=-1.298297475578389:1.3187811236619464,variable=\t]({1*0.2636871775680123*cos(\t r)+0*0.2636871775680123*sin(\t r)},{0*0.2636871775680123*cos(\t r)+1*0.2636871775680123*sin(\t r)});
\draw [line width=2pt,color=wrwrwr] (6,3)-- (7,3);
\draw [line width=2pt,color=wrwrwr] (7,3)-- (8,3);
\draw [shift={(10.38464384717023,2.487180641890687)},line width=2pt,color=wrwrwr]  plot[domain=0.694759916134297:1.8268778870334221,variable=\t]({1*0.8010286442300946*cos(\t r)+0*0.8010286442300946*sin(\t r)},{0*0.8010286442300946*cos(\t r)+1*0.8010286442300946*sin(\t r)});
\draw [shift={(10.395799263785687,3.5011563388713434)},line width=2pt,color=wrwrwr]  plot[domain=4.452523221142281:5.5907412170909865,variable=\t]({1*0.7849943984722737*cos(\t r)+0*0.7849943984722737*sin(\t r)},{0*0.7849943984722737*cos(\t r)+1*0.7849943984722737*sin(\t r)});
\draw [shift={(12.37223664729231,3.545636952422224)},line width=2pt,color=wrwrwr]  plot[domain=2.742787586115155:4.11370558442345,variable=\t]({1*0.6605148033435978*cos(\t r)+0*0.6605148033435978*sin(\t r)},{0*0.6605148033435978*cos(\t r)+1*0.6605148033435978*sin(\t r)});
\draw [shift={(11.575071066600975,3.5760175949446125)},line width=2pt,color=wrwrwr]  plot[domain=-0.935212957563345:0.314868569733304,variable=\t]({1*0.7157938726514839*cos(\t r)+0*0.7157938726514839*sin(\t r)},{0*0.7157938726514839*cos(\t r)+1*0.7157938726514839*sin(\t r)});
\draw [shift={(12.00915492892372,3.74885653575275)},line width=2pt,color=wrwrwr]  plot[domain=0.19557120262823546:2.9280039203807235,variable=\t]({1*0.2513102710947334*cos(\t r)+0*0.2513102710947334*sin(\t r)},{0*0.2513102710947334*cos(\t r)+1*0.2513102710947334*sin(\t r)});
\draw [shift={(13.571244562288568,2.5719892197896455)},line width=2pt,color=wrwrwr]  plot[domain=1.241332769475582:2.498565223945202,variable=\t]({1*0.7138021980356559*cos(\t r)+0*0.7138021980356559*sin(\t r)},{0*0.7138021980356559*cos(\t r)+1*0.7138021980356559*sin(\t r)});
\draw [shift={(13.597867360796997,3.467597203185726)},line width=2pt,color=wrwrwr]  plot[domain=3.805328598343443:4.99208466191838,variable=\t]({1*0.7590075925400742*cos(\t r)+0*0.7590075925400742*sin(\t r)},{0*0.7590075925400742*cos(\t r)+1*0.7590075925400742*sin(\t r)});
\draw [shift={(13.736432919397249,2.9920425106768644)},line width=2pt,color=wrwrwr]  plot[domain=-1.2982974755783845:1.3187811236619322,variable=\t]({1*0.2636871775680141*cos(\t r)+0*0.2636871775680141*sin(\t r)},{0*0.2636871775680141*cos(\t r)+1*0.2636871775680141*sin(\t r)});
\draw [line width=2pt,color=wrwrwr] (11,3)-- (12,3);
\draw [line width=2pt,color=wrwrwr] (12,3)-- (13,3);
\draw [shift={(10.273587194802213,3.004338428123557)},line width=2pt,color=wrwrwr]  plot[domain=1.9130770913733075:4.417617444885307,variable=\t]({1*0.2736215911040806*cos(\t r)+0*0.2736215911040806*sin(\t r)},{0*0.2736215911040806*cos(\t r)+1*0.2736215911040806*sin(\t r)});
\draw [shift={(5.273587194802215,3.0043384281235577)},line width=2pt,color=wrwrwr]  plot[domain=1.9130770913733142:4.417617444885297,variable=\t]({1*0.27362159110408074*cos(\t r)+0*0.27362159110408074*sin(\t r)},{0*0.27362159110408074*cos(\t r)+1*0.27362159110408074*sin(\t r)});
\draw [shift={(3.6948446751814736,1.623736272798756)},line width=2pt,color=wrwrwr]  plot[domain=2.6694198034150416:3.1132626430523826,variable=\t]({1*3.025936164257898*cos(\t r)+0*3.025936164257898*sin(\t r)},{0*3.025936164257898*cos(\t r)+1*3.025936164257898*sin(\t r)});
\draw [shift={(1.049821505780461,1.6984071960139773)},line width=2pt,color=wrwrwr]  plot[domain=-3.1706664948771905:0.011832889063124796,variable=\t]({1*0.3798593177773285*cos(\t r)+0*0.3798593177773285*sin(\t r)},{0*0.3798593177773285*cos(\t r)+1*0.3798593177773285*sin(\t r)});
\draw [shift={(-0.2792925281238718,1.8565347815086521)},line width=2pt,color=wrwrwr]  plot[domain=-0.08965813142236634:0.7293936897274242,variable=\t]({1*1.7158385933452578*cos(\t r)+0*1.7158385933452578*sin(\t r)},{0*1.7158385933452578*cos(\t r)+1*1.7158385933452578*sin(\t r)});
\draw [shift={(8.52915126122139,1.8622826911536126)},line width=2pt,color=wrwrwr]  plot[domain=2.5019305061325676:3.250289558167142,variable=\t]({1*1.9059654389688272*cos(\t r)+0*1.9059654389688272*sin(\t r)},{0*1.9059654389688272*cos(\t r)+1*1.9059654389688272*sin(\t r)});
\draw [shift={(5.243449632321056,1.7676839764164143)},line width=2pt,color=wrwrwr]  plot[domain=-0.0621222366083245:0.6117685708512434,variable=\t]({1*2.1457101333064292*cos(\t r)+0*2.1457101333064292*sin(\t r)},{0*2.1457101333064292*cos(\t r)+1*2.1457101333064292*sin(\t r)});
\draw [shift={(7.0109866130663745,1.6899037719532901)},line width=2pt,color=wrwrwr]  plot[domain=3.232657815689197:6.136060065160488,variable=\t]({1*0.3781191273019094*cos(\t r)+0*0.3781191273019094*sin(\t r)},{0*0.3781191273019094*cos(\t r)+1*0.3781191273019094*sin(\t r)});
\draw [shift={(10.710307706681078,1.728962955848445)},line width=2pt,color=wrwrwr]  plot[domain=-0.044131135139518385:0.5067599565066778,variable=\t]({1*2.6188214841202107*cos(\t r)+0*2.6188214841202107*sin(\t r)},{0*2.6188214841202107*cos(\t r)+1*2.6188214841202107*sin(\t r)});
\draw [shift={(14.352413141151859,1.8551246880328656)},line width=2pt,color=wrwrwr]  plot[domain=2.439109226539724:3.2784212255144256,variable=\t]({1*1.7719370147700173*cos(\t r)+0*1.7719370147700173*sin(\t r)},{0*1.7719370147700173*cos(\t r)+1*1.7719370147700173*sin(\t r)});
\draw [shift={(12.961808435034778,1.6651047952153624)},line width=2pt,color=wrwrwr]  plot[domain=3.282322847579891:6.142455113189487,variable=\t]({1*0.36841320986723747*cos(\t r)+0*0.36841320986723747*sin(\t r)},{0*0.36841320986723747*cos(\t r)+1*0.36841320986723747*sin(\t r)});
\draw (0.3894324115044247,3.3121956767200555) node[anchor=north west] {$i$};
\draw (1.7989124557522118,3.671470805941707) node[anchor=north west] {$j$};
\draw (3.1077153539823,3.2983774025192227) node[anchor=north west] {$k$};
\draw (0.7257000442477874,2.3034616600592646) node[anchor=north west] {$G'$};
\draw (1.6653534070796455,1.322364191800139) node[anchor=north west] {$G^i$};
\draw (5.456848761061946,3.3983774025192227) node[anchor=north west] {$i$};
\draw (6.7992107079646,3.6231068462387923) node[anchor=north west] {$j$};
\draw (8.141572654867256,3.305286539619639) node[anchor=north west] {$k$};
\draw (6.6992107079646,2.2896433858584317) node[anchor=north west] {$G'$};
\draw (10.457147013274335,3.243104305715892) node[anchor=north west] {$i$};
\draw (11.766627057522122,3.6816520236362943) node[anchor=north west] {$j$};
\draw (13.208989004424778,3.3121956767200555) node[anchor=north west] {$k$};
\draw (12.639162323008848,2.3449164826617626) node[anchor=north west] {$G'$};
\draw (6.766328805309733,1.2878185062980572) node[anchor=north west] {$G^j$};
\draw (11.833068008849557,1.37763728860347) node[anchor=north west] {$G^k$};
\begin{scriptsize}
\draw [fill=rvwvcq] (1,3) circle (2.5pt);
\draw [fill=rvwvcq] (2,3) circle (2.5pt);
\draw [fill=rvwvcq] (3,3) circle (2.5pt);
\draw [fill=rvwvcq] (6,3) circle (2.5pt);
\draw [fill=rvwvcq] (7,3) circle (2.5pt);
\draw [fill=rvwvcq] (8,3) circle (2.5pt);
\draw [fill=rvwvcq] (11,3) circle (2.5pt);
\draw [fill=rvwvcq] (12,3) circle (2.5pt);
\draw [fill=rvwvcq] (13,3) circle (2.5pt);
\end{scriptsize}
\end{tikzpicture}
 \caption{Graphs $G^i, G^j, G^k$ in \cref{cutcom}}
    \label{fig2}
\end{figure}

\bl \label{cutcom} Let $G$ be a tree, $G'$ be any connected graph and $G^v=G_v*G_w' $ with $v\in V(G), w\in V(G')$ and $\r^u(G^v)=\r (A_u)$ where  $A_u=\D(G^v)(u).$ If $i,j,k\in V(G)$ with $i \sim j\sim k$ then for any $u\in V(G)\cup V(G')$
\beq \r^u(G^i)+\r^u(G^k)\ge 2 \r^u(G^j). \eeq
Furthermore, either $\r^u(G^i)> \r^u(G^j)$ or $\r^u(G^k)> \r^u(G^j).$
\el
\pf Without loss of generality we have

\beq
\D(G^i)+\D(G^k)-2\, \D(G^j)=
\begin{blockarray}{ccccc}
G' & G_i & G_j & G_k \\
\begin{block}{(cccc)c}
0 & 0 & 2\, \J & 0 \\
0 & 0 & 0 & 0 \\
2\, \J^T & 0 & 0 & 0 \\
0 & 0 & 0 & 0 \\
\end{block}
\end{blockarray}
\eeq
where $G_i, G_j,G_k$ are components of $G-\{ij,\, jk \}$ containing $i,j,k$ respectively. 
Now since  $|V(G_j)|\ge 1$ and $|V(G')|\ge 2,$  for any $u\in V(G)\cup V(G')$ we get
\beq \D(G^i)(u)+\D(G^k)(u)-2\, \D(G^j)(u)\ge 0.
\eeq

Therefore if $\x$ is the Perron vector of $\D(G^j)(u),$ then
\beq \x^T \{ \D(G^i)(u)+\D(G^k)(u)-2\, \D(G^j)(u) \}\x \ge 0  \eeq
Thus
\be \r^u(G^i)+\r^u(G^k)-2 \, \r^u (G^j) &\ge \r(\D(G^i)(u)+\D(G^k)(u) )-2\, \r^u(G^j) \nonumber  \\
& \ge \x^T \{ \D(G^i)(u)+\D(G^k)(u)-2\, \D(G^j)(u) \}\x \nonumber  \\
\label{iden} &\ge 0. \qquad  [\text{ \cref{fig2} }]
\ee

Now if $|V(G_j)|\ge 2$ or if $u\ne j,$ then $\D(G^i)(u)+\D(G^k)(u)-2\, \D(G^j)(u)\ne 0$ and therefore Inequality \cref{iden} is strict and we are done. But if $u=j,$ then $\D(G^i)(u)+\D(G^k)(u)-2\, \D(G^j)(u)= 0.$ In this case to the contrary we assume that $\r^u(G^i)=\r^u(G^k)=\r^u (G^j)=\r \text{ ( say ).}$ Then we must have $\D(G^i)(u)\x =\r \x, \D(G^j)(u)\x =\r \x$ and $\D(G^k)(u)\x =\r \x.$ But we see that

\beq
\D(G^i)(u)-\D(G^j)(u)=
\begin{blockarray}{cccc}
G' & G_i  & G_k \\
\begin{block}{(ccc)c}
0 & -\J &  \J  \\
-\J^T & 0 & 0  \\
 \J^T & 0 & 0  \\
\end{block}
\end{blockarray}
\eeq

Therefore as $(\D(G^i)(u)-\D(G^j)(u))\x=0$ we observe from the above matrix that
\beq
{\ds \sum_{\ell \in V(G_i)}  }x_{\ell}=0 \text{ and } {\ds \sum_{\ell \in V(G_k)}  }x_{\ell}=0 \eeq
which are not possible as $\x$ is a positive vector and $V(G_i),\, V(G_k)$ are non empty. Hence either $\r^u(G^i)\ne \r^u(G^j)$ or $\r^u(G^k)\ne \r^u(G^j)$ and therefore we get the result as desired.
\qed

\bc \label{treeconvex} Let $T$ be a tree and $G'$ be a graph. Then $\r_2(G^i)$ is strictly quasiconvex on $T,$ where $G^i$ is as defined in \cref{cutcom}.
\ec
\pf Let $u$ be the vertex of $G^j,$ so that $\r_2(G^j)=\r(\D(G^j)(u)).$ Then by \cref{cutcom}, we get either $\r^u(G^i) >  \r_2(G^j) $ or $\r^u(G^k) > \r_2(G^j).$ But by \cref{r2def} we have $\r_2(G^i) \ge \r^u(G^i) $ and $\r_2(G^k)\ge \r^u(G^k).$
Hence the result follows.
\qed

As a result of repeated application of \cref{lem1} and \cref{treeconvex}  we get the following results.

\bt \label{r2mtpath} Among all trees of given order, the second largest distance Pareto eigenvalue is maximized in the path graph.
\et

\bt Among all trees of given order, the second largest distance Pareto eigenvalue is minimized in the star graph.
\et

From \cref{bound} we see that among all connected graphs of given order $n$ the minimum value of the second largest distance Pareto eigenvalue is $n-2$ and is uniquely achieved by the complete graph $K_n.$  Now using \cref{r2mtree,r2mtpath} we have the following result regarding the unique graph with maximum second largest distance Pareto eigenvalue among all connected graphs of given order.

\bt Among all connected graphs of given order, the second largest distance Pareto eigenvalue is maximized in the path graph.
\et

Let us reserve the following notations for the rest of this article.
For $S\subset [n],$ \beq X_n(S)&=\{ \x \in \R^n: \x \ge 0,\, \x^T\x=1, \,\x_i=0 \Leftrightarrow i\in S \}  \\
\text{ and }\quad X_n^k &=\{ \x \in X_n(S): S\subset [n], \, |S|=k \}. \eeq

\bt
If $A\in M_n$ is symmetric, non negative and irreducible (or positive), then
\beq
\r_2(A)={\ds \max_{\x\in X_n^1}\,\x^TA\x}.
\eeq
\et
\pf
From \cref{imp},
it can be observed that $\r_2(A)$ is the largest possible eigenvalue among all principal submatrices of $A$ of order $n-1.$
\par For any $i\in [n]$ let $B_i=A(i)$ and $\x \in \R^{n-1}$ be an arbitrary normalized vector. If we take $z\in \R^n$ such that $z_i=0$ and $z_j=\x_j$ for $j\neq i$, then
$z^TAz=\x^TB_i\x.$
Now
\beq
\r(B_i)&=\max_{\substack{\x\in \R^{n-1}\\
\x^T\x=1}}\x^TB_i\x \\
&=\max_{\x\in X_{n-1}^0}\x^TB_i\x \\
&=\max_{z\in X_n(\{i\})}z^TAz
\eeq

\beq
\text{ Therefore } \quad \r_2(A) &= \max_{i\in [n]}\r(B_i) \\
&= \max_{i\in [n]} \max_{z\in X_n (\{ i \}) }z^TAz \\
&= \max_ {\substack{z\in X_n (\{ i \}) \\ i\in [n] }}    z^TAz \\
&= \max_{z\in X_n^1} z^TAz
\eeq
This completes the proof.
\qed

\bc If $z\in X_n^1$ and   $A\in M_n$ is symmetric, non negative irreducible (or positive) then \beq  \min _{z\in X_n^1}z^TAz \le  \r_2(A)\le \max_{z\in X_n^1}   z^TAz,\eeq  with left hand (right hand) equality if and only if $z$ is the normalized Pareto eigenvector of $A$ corresponding to $\r_2(A).$ \ec

\bc If $A\in M_n$ is symmetric, non negative irreducible (or positive) then  \beq \r_2(A)=\max_{x\, \in T_n^1}\frac{x^T A x}{x^Tx}  \eeq where $T_n^k=\{x\in \R^n: x\ge 0 \mbox{ and } x_i=0 \mbox{ for exactly k values of } i\in [n]\}.$  \ec

\bt If $G$ is a connected graph of diameter $d$ and $T_{min}={\ds \min_{v\in V(G)}}Tr(v),$ then
\beq \r_2(G)\ge \frac{T_{min}-2d+\sqrt{(T_{min}-2d)^2+4(n-d-1) } }{2} \eeq
Equality holds if and only if $G=K_n.$
\et
\pf Let $\x=( x_1, x_2, \ldots , x_n)^T $ be the Pareto eigenvector of $\D(G)$ corresponding to $\r_2(G)$ with $x_i=0$ i.e. $\r_2(G)=\r(\D(G)(i)).$ \\
Let \beq x_j=\min_{\ell \ne i}x_{\ell} \hs \text{ and } x_k= \min_{\ell \ne i,j}x_{\ell}.\eeq

For $j\in V(G)$ if $T_j=Tr(j),$ then from Pareto eigenequations we have
\be \r_2 x_j &= {\ds \sum_{\ell \in V(G)}d_{jl}x_{\ell}} \nonumber \\
\label{Huv1} & \ge (T_j-d_{ji})x_k \\
\label{Huv2} & \ge (T_j-d)x_k
\ee

\be \r_2 x_k &={\ds \sum_{\ell \in V(G)}} d_{k\ell}x_{\ell} \nonumber \\
& \ge d_{kj}x_j +(T_k-d_{ki}-d_{kj})x_k \nonumber \\
\label{Huv3} & \ge x_j +(T_k-2d)x_k
\ee

From \cref{Huv2,Huv3} we get
\be \r_2 (\r_2-T_k+2d)-T_j+d & \ge 0\nonumber \\
\Rightarrow \quad \r_2 & \ge \frac{T_k-2d+ \sqrt{(T_k-2d)^2+4(T_j-d)} }{2}\nonumber \\
\label{Huv4} & \ge \frac{T_k-2d+\sqrt{(T_k-2d)^2+4(n-d-1) }}{2}.
\ee
\par If $G=K_n,$ then $T_k=n-1$ and $d=1.$ Therefore by \cref{k} equality holds in \cref{Huv4}.
 \par Now suppose the equality holds in \cref{Huv4} then equality must hold in all the inequalities \crefrange{Huv1}{Huv4}.
  Equality in \cref{Huv4} gives $d_j=n-1.$ Again equality in \cref{Huv3} implies $d_{jk}=d.$ Therefore we must have $d=1$ i.e. $G=K_n.$
\qed

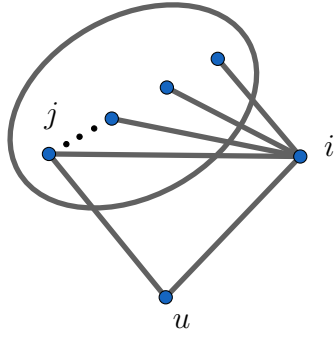
\begin{figure}
    \centering
\begin{tikzpicture}[line cap=round,line join=round,>=triangle 45,x=1cm,y=1cm]
\clip(1.5, 1) rectangle (6.53232323232323,5.7);
\draw [rotate around={29.419076151748374:(3.491807440805782,4.1009095965828175)},line width=2pt,color=wrwrwr] (3.491807440805782,4.1009095965828175) ellipse (1.7564440375567605cm and 1.1938784256199682cm);
\draw [line width=2pt,color=wrwrwr] (3.9221191003533216,1.5612270569002777)-- (2.369622132573963,3.468098332542317);
\draw [line width=2pt,color=wrwrwr] (3.9221191003533216,1.5612270569002777)-- (5.710865606707801,3.4343483984601573);
\draw [line width=2pt,color=wrwrwr] (5.710865606707801,3.4343483984601573)-- (4.613992749037602,4.7337208606233165);
\draw [line width=2pt,color=wrwrwr] (5.710865606707801,3.4343483984601573)-- (3.9374277369201427,4.35219923724881);
\draw [line width=2pt,color=wrwrwr] (3.20473901932965,3.9390289077804863)-- (5.710865606707801,3.4343483984601573);
\draw [line width=2pt,color=wrwrwr] (5.710865606707801,3.4343483984601573)-- (2.369622132573963,3.468098332542317);
\draw (3.8737373737373733,1.4832402234636903) node[anchor=north west] {$u$};
\draw (5.885521885521885,3.864525139664808) node[anchor=north west] {$i$};
\draw (2.1649831649831643,4.2946927374301715) node[anchor=north west] {$j$};
\begin{scriptsize}
\draw [fill=rvwvcq] (2.369622132573963,3.468098332542317) circle (2.5pt);
\draw [fill=rvwvcq] (4.613992749037602,4.7337208606233165) circle (2.5pt);
\draw [fill=rvwvcq] (5.710865606707801,3.4343483984601573) circle (2.5pt);
\draw [fill=rvwvcq] (3.9221191003533216,1.5612270569002777) circle (2.5pt);
\draw [fill=rvwvcq] (3.9374277369201427,4.35219923724881) circle (2.5pt);
\draw [fill=rvwvcq] (3.20473901932965,3.9390289077804863) circle (2.5pt);
\draw [fill=black] (2.6010481839584783,3.598601744977194) circle (1pt);
\draw [fill=black] (2.7763505298699602,3.6974564513182555) circle (1pt);
\draw [fill=black] (2.9710025581143404,3.8072226326590712) circle (1pt);
\end{scriptsize}
\end{tikzpicture}
    \caption{Graph $H_{i,u}$ as in \cref{last} }
    \label{fig4}
\end{figure}

\bt \label{last} Let $G$ is a connected graph and $(\r_2(G),\x)$ be a distance Pareto eigenpair of $G.$ If the second largest component of $\x$ corresponds to $j\in V(G)$ then
\beq \r_2(G)\le  \frac{T_{j}-2+\sqrt{(T_{j}-2)^2+4(n-2) } }{2}. \eeq
If equality holds, then $G \cong H_{u,v}$ for some $u,v \in V(G).$
\et
\pf Let $\x=( x_1, x_2, \ldots , x_n)^T $ be the Pareto eigenvector of $\D(G)$ corresponding to $\r_2(G)$ with $x_u=0$ i.e. $\r_2(G)=\r(\D(G)(u)).$ \\
Let \beq x_i=\max_{k\in [n]}x_k \hs \text{ and } x_j= \max_{k\ne i}x_k \eeq
From Pareto eigenequations we have
\be
\r_2x_i & = {\ds \sum_{\ell \in V(G) } }d_{i\ell}x_{\ell} \nonumber \\
  &\le (T_i-d_{iu})x_j \nonumber \\
\label{r2max1} & \le (T_i-1)x_j
\ee
and
\be \r_2 x_j &\le d_{ij}x_i+(T_j-d_{ju}-d_{ij} )x_j \nonumber \\
\label{r2max2} & \le x_i +(T_j-2)x_j \ee

Now \cref{r2max1,r2max2} together implies
 $\r_2^2-(T_j-2)\r_2-(T_i-1) \le 0.$ Which in turn gives
\be \label{r2max3} \r_2 \le \frac{T_j-2+\sqrt{(T_j-2)^2+4(n-2) }}{2}
\ee
Now suppose the equality holds in \cref{r2max3},
then equality must hold in \cref{r2max1,r2max2} as well. \\

Equality in \cref{r2max1} implies
\be   i \sim u \text{ and } x_k=x_j \hs \text { for all } k\ne i,u \nonumber \\
 \label{r2max5} \text{ Thus }  T_k=T_j \hs \text { for all } \, k\ne i,u \qquad [ \text{ using \cref{r2max2} }]
\ee

Equality in \cref{r2max2} implies $i\sim j\sim u $ and equality in \cref{r2max3} gives $d_i=n-1$ and so diam$(G) \le 2.$ Hence using \cref{r2max5} we get  $T_k =d_k+2(n-d_k-1)
=2(n-1)-d_k$ for all  $k\ne i,u.$

Combining all the above arguments, we get the required result.
\qed

\bigskip



\bibliographystyle{abbrv}


\end{document}